\documentclass[12pt]{amsart}
\usepackage{amsmath,amsfonts,amsthm,amstext,amscd}
\usepackage[dvips]{graphicx}
\usepackage{array}
\usepackage{hyperref}
\usepackage{newlfont}

\usepackage[all]{xy}

\newcommand{\spand}{Span\{A(t),\dot{A}(t),...,A(t)^{(k-2)}\}}

\newcommand{\jusk}{( A(t)|\dot{A}(t)|...|A^{(k-1)}(t))}
\newcommand{\dem}{\textit{Proof:}}

\newcommand{\obs}{\textit{Remark.}}

\newcommand{\lieglkn}{\mathfrak{gl}(kn)}

\newtheorem{The}{Theorem}[section]

\newtheorem{lem}[The]{Lemma}
\newtheorem{Prop}[The]{Proposition}
\newtheorem{Def}[The]{Definition}

 \numberwithin{equation}{section}

\begin{document}

%
%
%
%
%
%
%
%
%

\title[Geometry of fanning curves]{Geometry of fanning curves in divisible Grassmannians }

\author[Dur\'an]{Carlos E. Dur\'an}

\address{C. Dur\'an\\Departamento de Matem\'atica, UFPR \\
 Setor de Ci\^encias Exatas, Centro Polit\'ecnico, \\
 Caixa Postal 019081,  CEP 81531-990, \\
 Curitiba, PR, Brazil}

\email{cduran@ufpr.br}

\author[Peixoto]{C\'{i}ntia R. de A. Peixoto}

\address{C.R.A. Peixoto, IMECC-UNICAMP, Pra\c{c}a Sergio Buarque de Holanda, 651, 
Cidade Universit\'aria - Bar\~ao Geraldo, 
Caixa Postal: 6065 
13083-859 Campinas, SP, Brasil }

\curraddr{Laboratoire Paul Painlev\'e Bat. M3, Universit\'e
des Sciences et Technologies, 59 655 Villeneuve d’Ascq, France}

\email{cintiarap@yahoo.com.br}

\thanks{C.R.A. Peixoto was financially supported by CNPq, grant \# 234165/2014-6}



\date{}

\maketitle
\begin{abstract}

We study the geometry of fanning curves in the Grassmann manifold of n-dimensional 
subspaces of $\mathbb{R}^{kn}$; we construct a complete system of invariants 
which solve the congruence problem. The geometry of the invariants themselves 
and their relation with classical invariants is also studied.  
\end{abstract}

\section{Introduction}

Consider the {\em divisible Grassmanians} $Gr(n,kn)$ of $n$-dimensional subspaces of $\mathbb{R}^{kn}$ 
as a homogenous space of $GL(kn)$. The aim of this paper is to study the geometry of curves 
$\ell(t) \in Gr(n,kn)$; these curves representes the projectivized geometry of solutions of 
systems of $n$ linear ordinary differential equations of order $k$. 
We will construct and study complete invariants that solve the congruence 
problem; but the main thrust of this paper is a thorough investigation of 
the equivariant geometry of the spaces of jets of curves in the divisible 
Grassmannian, by modelling them as adjoint orbits in the Lie algebra $\lieglkn$. 
Both the invariants and their geometric interpretation are a consequence of the 
adjoint model. 

This work extends several Klein geometries:

\begin{itemize}  
 \item The classical projective invariants of ordinary differential equations studied 
by Wilczynski (\cite{wil}), an important distinction between our invariants and the 
Wilczynski invariants
 is that he considers 
a single differential equation whereas we consider systems; this is reflected in the {\em non-commutativity}
of the invariants. 
\item Our moving frames generalize the ``commutative'' case $n=1$, that is, 
the linear geometry of curves in the real projective 
space, studied by Cartan (\cite{cartan}).
\item The main inspiration is the paper \cite{duran} that studies the case $k=2$ of systems of 
second order linear differential equations. In the general case treated here, in addition to extra combinatorial complexity  some new phenomena appear; for example,
 the natural matrix of invariants of proposition \ref{Jacobi-in-base}
 does {\em not} coincide with the pullback of the Maurer-Cartan form. 
\item The common denominator of all these cases is the work of Flanders (\cite{flanders}) of curves in 
$\mathbb{R}P^1$. 
\end{itemize}

The case $k=2$ (the ``half-Grassmannian'') was studied in \cite{duran}; we 
briefly describe that paper here:  the main insight of \cite{duran} is that the linear invariants of curves in the Grassmannian and their geometry are completely
 described by the {\em fundamental endomorphism} $F$ and its derivatives, which is an equivariant 
map
 from $1$-jets of curves in the Grassmannian into the Lie algebra $\mathfrak{gl}(2n)$ endowed with 
the Adjoint action. The first derivative of the fundamental endomorphism is a reflection whose 
$-1$-eigenspace is the curve $\ell(t)$ itself, and the $1$-eigenspace furnishes an
equivariant complement $h(t)$ of $\ell(t)$ which is called the {\em horizontal curve}. 
The main geometric invariant of a fanning curve $\ell(t)$, described in \cite{duran}, 
is its Jacobi endomorphism, that describes how the horizontal curve moves with respect to 
$\ell(t)$, and it gives the natural notion of curvature for fanning curves in the 
Grassmannian $Gr(n,\mathbb{R}^{2n})$. There is a close relationship between the matrix 
generalization of the Schwarzian derivative  (based on the work of  Zelikin \cite{Zelikin}) and the Jacobi endomorphism, also studied in \cite{duran}.

Following the same spirit as the half-Grassmannian case, our study will proceed along the following  main lines:

\medskip

\noindent{\bf Fanning frames and fanning curves:} We study curves $\ell(t) \in Gr(n,\mathbb{R}^{kn})$ via frames $A(t)$ spanning $\ell(t)$:

\begin{Def}
A frame $A(t)$, organized as a curve of $kn \times n$ matrices is fanning,
if the $kn \times kn$ matrix $\mathbf{A}(t) := ( A(t)|\dot{A}(t)|...|A^{(k-1)}(t))$ formed juxtaposing $A(t)$ and its derivatives is invertible  for all $t$.  This condition
depends only on the space $\ell(t)$ spanned by the columns of $A$; thus we say that a curve
$\ell(t)$ is fanning if a curve of frames $A(t)$ spanning $\ell(t)$ is fanning.
\end{Def}

Fanning curves form an open and dense subset of all differentiable curves, and therefore
it is a natural non-degeneracy condition. For the rest of this work, we will always assume to work in the 
set of fanning curves. Observe that the matrix $\mathbf{A}(t)$ gives a (highly non-canonical) $GL(kn)$ equivariant 
lift of the curve $\ell(t)$ into $GL(kn)$, the later endowed with the canonical left action on itself.  

Another construction defined for (fanning) frames that only depends on the curve is the canonical flag 
\[
Span\{A(t)\} \subset Span\{A(t),\dot{A}(t)\}\subset  \dots
\subset Span\{A(t),\dot{A}(t),...,A(t)^{(k-2)}\} 
\]
$$Span\{A(t),\dot{A}(t),...,A(t)^{(k-2)}\} \subset \mathbb{R}^{kn} \, .$$

We call the last non-trivial space $Span\{A(t),\dot{A}(t),...,A(t)^{(k-2)}\}$ of this sequence the 
{\em vertical space} $v(t)$.

 The next items correspond to the sections of this paper:

\medskip

\noindent{\bf Normal forms:} We will construct a 
{\em normal form} for frames spanning the given curve $\ell(t)$. This normal form gives a canonical 
way of extending an initial frame $A(0)$ of $\ell(0)$ to a frame $A(t)$ spanning $\ell(t)$; linear 
relations between the derivatives of a normal frame furnish a complete system of  invariants, which generalizes  for systems the Wilczynski 
invariants for single differential equations ( \ref{mainCongruence}). However, for systems the invariants are matrices, instead of numbers, and the non-commutativity implies
that there is an `up to conjugation by a constant matrix" in the conjugacy theorem. Therefore, the actual invariants are the linear transformations expressed as matrices in a given basis. If actual matrix invariants are wanted, it is necessary to further 
normalize the curve.
\medskip

\noindent{\bf The fundamental endomorphism and its derivatives:}  We 
generalize the fundamental 
endomorphism of \cite{duran}, obtaining an Adjoint-equivariant map $F(t)$ into $\lieglkn$. 
The ``derivative'' 
\[
D(t)= \frac{1}{k} (2 \dot{F}(t) - (k-2)I)
\]
 is the {\em fundamental reflection}, 
whose $-1$ eigenspace 
is the vertical $v(t)$; its $+1$ eigenspace, the {\em horizontal curve 
$h(t)$}, will be 
a fundamental piece of the study of the invariants, together with 
the related {\em horizontal derivative} $H(t)$ which spans $h(t)$. 
The horizontal derivative has the form 
\[
H(t)= A(t)^{k-1} + \text{extra terms depending on lower order derivatives of } A .
\]

Thus the horizontal derivative has the same spirit as a $k-1$-th order ``covariant" derivative. 
An important remark is that, for $k=2$, the horizontal derivative in a normal frame is just 
the ordinary derivative of the frame, whereas this fails for $k>2$. This influences, for example the  
Cartan lift of the $k-1$-jet of a curve to $GL(kn)$: there are two choices, 
one of them using 
$( A(t),\dot{A}(t),...,A^{(k-1)}(t))$ for a normal frame, and the other, still 
with a normal frame but
using the horizontal derivative,  $( A(t),\dot{A}(t),...,A(t)^{(k-2)}, H(t))$ (these two lifts coincide 
for $k=2$). 

Taking one more derivative, we arrive at the 
matrix of invariants, the {\em Jacobi endomorphism}, whose entries are in direct relationship with the normal form invariants and has the geometric interpretation of 
measuring the velocity of canonically defined curves of flags associated to the 
curve in the Grassmannian. 

\medskip

\noindent{\bf Geometry of jets of fanning curves in the Grassmannian:} Here we see how the 
invariants arise naturally by representing the prolonged action of $GL(kn)$ on jets of curves on 
the Grassmannian as the Adjoint action on $\lieglkn$. We shall see that once one is comitted to 
an Adjoint representation, there is essentially no choice of invariants. 

Also, we use this Adjoint representation to give a better understanding of 
the spaces of jets of curves in the divisible Grassmannian as a  $GL(kn)$-space.
In particular, by restricting to ``standard" curves, we give numerical invariants
that serves as coordinates for the space of orbits on the $k+1$-jets (the first level on which that $GL(kn)$-action fails to be transitive).

\section{Normal Frames} \label{normal}

The fanning condition for frames $A(t)$ spanning curves
$\ell(t) \in Gr(n,\mathbb{R}^{kn})$, means that at each instant $t$ the columns of
$A(t), \dot{A}(t), \ddot{A}(t), \cdots, A(t)^{(k-1)}$ span $\mathbb{R}^{kn}$ and therefore we
can write  $A^{(k)}$ as a linear combination of $A(t), \dot{A}(t), \ddot{A}(t)\cdots A(t)^{(k-1)}$.
This gives a system of order $k$ differential equations satisfied by the columns of $A$;
in this part we will adopt the notation of  Wilczynski \cite{wil} for writing the
coefficients:
for example, in the case $k=3$ we write
\begin{equation}\label{eq3}
A^{(3)} + 3\ddot{A}P_1(t) + 3\dot{A}P_2(t) + AP_3(t) = 0,
\end{equation}
where $P_1(t)$, $P_2(t)$ and $P_3(t)$ are smooth curves of $n \times n$ matrices,
and in the general case $Gr(n,\mathbb{R}^{kn})$, we have
\begin{equation}\label{eqk}
A^{(k)} + {k \choose 1}A^{(k-1)}P_1 + {k \choose 2}A^{(k-2)}P_2 +...+{k \choose k-1}\dot{A}P_{k-1}+ AP_{k} = 0,
\end{equation}
where $P_i$ for all $i$ are smooth curves depending on $t$ of $n \times n$ matrices.

\medskip

\noindent \textit{Remark.} This is the first instance of giving the case $n=3$ first and then the general
case.  It is much easier to visualize the combinatorics in this case, but most of that, 
this is the first $n$ where the differences from  $Gr(n,\mathbb{R}^{2n})$ (\cite{duran}) appears. 

\medskip

\begin{Def}
A  frame $A(t)$ of a fanning  curve in $Gr(n,\mathbb{R}^{kn})$ is said to be normal if the columns of its kth-derivative $A^{(k)}$ are linear combinations of the columns of all derivatives of order equal or less than $k-2$, for all values of $t$.
\end{Def}

This definition is motivated by the normal frames defined by Cartan \cite{cartan}, and coincides
for $k=2$ with the normal frames defined in \cite{duran}.

In general, when two frames  $A(t), B(t) $  span a curve $\ell(t) \in Gr(n,\mathbb{R}^{kn})$  this means that there is  $X(t) \in Gl(n)$ that  $A(t)=B(t)X(t)$. 
If $\ell(t)$ is a fanning curve of $n$-dimensional subspaces in $\mathbb{R}^{kn}$, there is a normal frame that spans it. In order to obtain a normal frame for $\ell(t)$ we use a method of reduction for differential equations of order $k$, described for Wilczynski \cite{wil}, that consists in a change of variables that results in a equation without the term of order $k-1$.

If in equation \ref{eqk} we put $A(t)=B(t)X(t)$, where $X(t)$ satisfies $\dot{X}(t)=-X(t)P_1(t)$ with $X(0)=I$, we obtain the equation
$$B^{(k)} + {k \choose 2}B^{(k-2)}Q_2(t) +...+{k \choose k-1}\dot{B}Q_{k-1}(t)+ BQ_{k}(t) = 0,$$
where :
$$Q_j(t)= \sum_{i=0}^j {j \choose i} \left( \frac{d^{j-i}}{dt^{j-i}}X \right)P_iX^{-1}$$

with $j=2,...,k$ and $P_0=1.$

We find in particular:
\begin{eqnarray*}
Q_2 &= &X (P_2 - P^2_1 -\dot{P}_1)X^{-1},\\
Q_3& = & X(P_3 - 3 P_1P_2 -2\dot{P}_1P_1+2P_1\dot{P}_1 +2P^3_1 -\ddot{P}_1)X^{-1},\\
Q_4 & = & X ( P_4 -4P_1P_3 +6P^2_1P_2 -6\dot{P}_1P_2 +3\dot{P}_1P^2_1-3P^2_1\dot{P}_1+{} \nonumber\\
& & {} +6P_1\dot{P}_1P_1+3P_1\ddot{P}_1 -3\ddot{P}_1P_1 -3P^4_1 +3\dot{P}^2 -P^{(3)})X^{-1}.
\end{eqnarray*}

This matrices $Q_j$ in $Gl(n)$ are pointwise conjugate under the action of $Gl(n)$ in the space of frames, that is  $Q_j = X T_j X^{-1}.$ (In the case $n=1$ these are the semi-invariants of \cite{wil}.)

If $A(t)$ is a fanning frame that satisfies equation (\ref{eqk}), let us define the $Schwarzian$ of $A(t)$ as the function
$$\{A(t),t\} = 2(P_2(t) - P_1(t)^2 -\dot{P}_1(t)).$$


\medskip

\noindent\textit{Remark.} The notation adopted here does not change the results in \cite{duran},  the case $k=2$. When the fanning frame is of the form $A(t)= \left( \begin{array}{c}
I \\
M(t)
\end{array}\right),$ we still have  $\{A(t),t\}= \frac{d}{dt}(\dot{M}^{-1}\ddot{M})-(1/2)(\dot{M}^{-1}\ddot{M})^2.$  And the normal form of the equation \ref{eqk} still be
$$\ddot{A}+(1/2)A\{A(t),t\}=0.$$ 
The only change is the form of the horizontal derivative that changes to $ H(t) = \dot{A}(t) + A(t)P_1(t)$, without affecting the other results. In this way, this form generalizes  the case $Gr(n,\mathbb{R}^{2n})$, and in this work the Schwarzian is the first invariant.
We will use the notation $h_{j-2}[A(t),t]$ or simply $h_{j-2}$ to denote $X^{-1}(t)Q_j(t)X(t)$, with $j \in \{3,...,k\}$. For example,
\[
h_1[A(t),t] = X^{-1}(t)Q_2(t)X(t) = P_3 - 3 P_1P_2 -2\dot{P}_1P_1+2P_1\dot{P}_1 +2P^3_1 -\ddot{P}_1 \, .
\]

We emphasize the Schwarzian by calling it $\kappa$ instead of $h_0$. 

The following properties of the Schwarzian and of the coefficients $h_j$ follow from the reduction of the equation:

\begin{Prop}\label{actschw}
Let $A(t)$ be a fanning frame.

(1) If $X(t)$ is a smooth curve on $Gl(n)$, then 
\[\{A(t)X(t),t\} = X(t)^{-1}\{A(t),t\}X(t)\]  and $h_{j}[A(t)X(t),t] =  X(t)^{-1}h_{j}[A(t),t]X(t)$, for $ j \in \{1,...,k-2\}$.

(2) If $T$ is a transformation on $Gl(kn)$, then $\{TA(t),t\} = \{A(t),t\}$ and $h_{j}[TA(t),t] =  h_{j}[A(t),t]$, for $ j \in \{1,...,k-2\}$.
\end{Prop}

\begin{Prop}\label{cteX}
Let $\ell(t)$ be a fanning curve of $n$-dimensional subspaces in $\mathbb{R}^{kn}$.  If $A(t)$ e $B(t)$  are two normal frames spanning $\ell(t)$, there is  a fixed invertible $n\times n$ matrix $X$ such that $B(t)= A(t)X$.
\end{Prop}

\dem

If $A(t)$ e $B(t)$ are two normal frames spanning $\ell(t)$ in $Gr(n,\mathbb{R}^{kn})$, then there is a curve of $n\times n$ invertible matrices $X$ such that $B(t)= A(t)X(t).$  Differentiating the equation $k$-times

$$B^{(k)} = A^{(k)}X +{k \choose 1}A^{(k-1)}\dot{X} + ...+{k \choose k-1}\dot{A}X^{(k-1)}+ AX^{(k)}.$$

Observe that $A(t)$ is normal, then $A^{(k)}$ depends only on the $k-2$-derivatives of $A(t)$; but $B(t)$ is normal too, so the only possible way in which the columns of $B^{(k)}$ could be linear combinations of the columns of the $k-2$-derivatives 
of $B(t)$ is that $\dot{X}$ be identically zero.  \qed

\medskip

Proposition \ref{cteX} has two important consequences: first, the juxtaposed matrix $\mathbf{A}(t)$ in the 
introduction is almost canonically defined for normal frames: it just depend on the choice of a basis 
of $\ell(0)$.

Also, as mentioned in the introduction, a fanning curve 
$\ell(t) \in Gr(n,kn)$ naturally produces a linear flag 
\[
Span\{A(t)\} \subset Span\{A(t),\dot{A}(t)\}\subset  \dots \subset   Span\{A(t),\dot{A}(t),...,A(t)^{(k-2)}\} 
\]
\[
Span\{A(t),\dot{A}(t),...,A(t)^{(k-2)}\}  \subset \mathbb{R}^{kn} \, ,
\]
but now proposition \ref{cteX}   makes this flag a {\em decomposition} 
flag:
\[
Span\{A(t)\}\oplus Span\{\dot{A}(t)\}\oplus \dots \oplus Span\{A(t)^{(k-1)}\}
\]
In general, transforming a linear flag of nested subspaces into a decomposition is only possible in the 
presence of an Euclidean structure by taking complements, but here the normal frame construction on fanning curves gives 
this additional structure. 

We shall see, however, that neither the lift $\mathbf{A}(t)$ nor the decomposition are  the most convenient ones; 
the ``right'' constructions
will be given in sections \ref{horcur} and \ref{horder} by means of the horizontal curve and the horizontal derivative.

Let us now prove the main result of this section, which essentially solves the congruence problem:
\begin{The} \label{mainCongruence}
Two fanning curves of $n$-dimensional subspaces of $\mathbb{R}^{kn}$ are congruent if and only if there exists a constant 
$n\times n$ invertible matrix $X$ such that the Schwarzians and the matrices $h_j$, for $j=1,..,k-2$, of any two of their normal frames are conjugate by $X$. 
\end{The}

\dem

Let $A(t)$ and $B(t)$ be two normal frames spanning congruent fanning curves. Then there is a linear transformation $T$ of $\mathbb{R}^{kn}$ such that $TA(t)$ and $B(t)$ span the same curve. Since $TA(t)$ is a normal frame too, using the Proposition \ref{cteX}, there exists $X$ constant such that $TA(t)=B(t)X$. And Proposition \ref{actschw} tell us that the Schwarzian and the $h_j$, for all $j$, of $A(t)$ and $B(t)$ are conjugate by a constant matrix on $Gl(n)$.

On the other hand, let $A(t)$ and $B(t)$ be normal frames such that $\{B(t),t\} = X(t)^{-1}\{A(t),t\}X(t)$ and $h_{j}[B(t),t] =  X(t)^{-1}h_{j}[A(t),t]X(t)$, for all $ j \in \{1,...,k-2\}$, then we can consider, without loss of generality, that $A(t)$ and $B(t)$ are two normal frames with the same Schwarzian and $h_j$, for all $j$. Assume that 
\[
T =( A(0)|\dot{A}(0)|...|A^{(k-1)}(0))( A(0)|\dot{A}(0)|...|A^{(k-1)}(0))^{-1};
\]
since $A(t)$ is  normal, then $D(t)= TA(t)$ satisfies a differential equation with order $k$:
$$ D^{(k)} + 1/2 D^{(k-2)}\{B(t),t\}+ D^{(k-3)}h_{1}[B(t),t]+ ...+ D(t) h_{k-2}[B(t),t]= 0.$$
Therefore $D(t)$ and $B(t)$ satisfy the same differential equation of order $k$ and with the same initial conditions. It follows that $D(t)=B(t)$ and then $A(t)$ is congruent to $B(t) $. \qed

\section{The fundamental endomorphism and its derivatives} \label{fundamentaletc}

\subsection{The fundamental endomorphism}

\begin{Def}
The fundamental endomorphism of a fanning frame $A(t)$ at a given time $t$ is the linear transformation
$\mathbb{R}^{kn} \rightarrow \mathbb{R}^{kn}$  defined by the equations $F(t)A(t)= 0$, $F(t)\dot{A}(t)=A(t)$, 
$F(t)A^{(2)}(t)= 2 \dot{A}(t)$, ..., $F(t)A^{(k-1)}(t)= (k-1) A^{(k-2)}(t)$.

\end{Def}

\medskip

\noindent{\bf Remark.}
Equivalently, we could have defined the fundamental endomorphism by the transformation, 
defined in the canonical basis, associated to the matrix $F$ of the theorem \ref{endfund} below.


\medskip
The fundamental endomorphism does not depend on the fanning frame (that is  1 of \ref{actF}  below), 
therefore it is defined for fanning curves in the Grassmannian. Furthermore, if $\ell(t)$ is a fanning curve spanned by $A(t)$ and $F(t)$
its fundamental endomorphism, then $$\ell(t)=Span\{A(t)\}= Im \{F(t)^{k-1}\} \subset Span\{A(t),\dot{A}(t)\}=Im\{F(t)^{k-2}\},$$
and we have that, for all
$i \in \{1,...,k-2\}$ $$ Span\{A(t),...,A(t)^{(i)}\}= Im \{F(t)^{k-i}\} \subset Span\{A(t),...,A(t)^{(i+1)}\},$$
$$\text{and} \quad  Span\{A(t),...,A(t)^{(i+1)}\}=Im\{F(t)^{k-(i+2)}\}$$
moreover $Span\{A(t),...,A(t)^{(i)}\}$ does not depend on the frame, for all $i$. In the case of normal frames, $Span\{A(t)^{(i)}\}$,
for all $i$, does not depend on the frame too.

\begin{Prop}\label{actF}
Let $A(t)$ be a fanning frame. Its fundamental endomorphism $F(t)$  satisfies the following properties:
\begin{enumerate}
\item If $X(t)$ is a smooth curve on  $Gl(n)$, the fundamental endomorphism of $A(t)X(t)$ is $F(t)$.
\item If $T$ is a matrix on $Gl(kn)$, the fundamental endomorphism of $TA(t)$ is $TF(t)T^{-1}$.

\end{enumerate}
\end{Prop}

\dem

The proof is the same as the half-grassmannian $Gr(n,2n)$ of \cite{duran}. \qed

\subsection{The fundamental reflection and the horizontal curve}\label{horcur}
We now take derivatives of the fundamental endomorphism and study the resultant geometry. 

\begin{Prop}\label{d}
Let $F(t)$ be the fundamental endomorphism of a fanning frame $A(t)$. At each value of 
$t$, $D(t)= \frac{1}{k} (2 \dot{F}(t) - (k-2)I)$ is a reflection whose $-1$ 
eigenspace is the vertical space $v(t)$.
\end{Prop}

\dem

We first observe that differentiating the identities 

$F(t)A(t)=0, F(t)\dot{A}(t)=A(t),..., F(t)A(t)^{(k-2)}=(k-2)A(t)^{(k-3)}$, we obtain, 
respectively, that 
$$\dot{F}(t)A(t) = -A(t), \dot{F}(t)\dot{A}(t) = -\dot{A}(t),..., \dot{F}(t)A(t)^{(k-2)} = -A(t)^{(k-2)}.$$ 
Consequently $$D(t)A(t)=-A(t), D(t)\dot{A}(t)=-\dot{A}(t), ..., D(t)A(t)^{(k-2)}=-A(t)^{(k-2)}.$$ 
Since the range of $F(t)$ is $Span\{A(t), \dot{A}(t),...,  A(t)^{(k-2)}\}$, then $\dot{F}(t)F(t)=-F(t)$.

Now we show that  $D(t)^2=I$. Differentiating $F(t)A(t)^{(k-1)}=(k-1)A(t)^{(k-2)}$ and using 
that $\dot{F}(t)F(t)=-F(t)$, we have:
$$\dot{F}(t)^2A(t)^{(k-1)} - F(t)A(t)^{(k)} = (k-1)\dot{F}(t)A(t)^{(k-1)},$$

then
$$(\dot{F}(t)^2 -(k-2)\dot{F}(t))A(t)^{(k-1)} = (F(t)A(t)^{(k-1)})', $$
and consequently,
$$(\dot{F}(t)^2 -(k-2)\dot{F}(t))A(t)^{(k-1)} = (k-1)A(t)^{(k-1)}.$$

Now multiplying by 4, using that $ 4k-4 = k^2 - (k-2)^2, \forall t$ and  completing the square we obtain that
$$\frac{1}{k^2} (2 \dot{F}(t) - (k-2)I)^2 A(t)^{(k-1)} = A(t)^{(k-1)}.$$ \qed

\medskip

It is  useful to think in terms of the fundamental projection $P(t) := \frac{(I-D(t))}{2}$ associated to a 
the fundamental reflection; $P(t)$ has the vertical space $v(t)$ as its image.  Its kernel is distinguished:

\medskip 

\begin{Def}
Let $\ell(t)$ be a fanning curve and let $F(t)$ be its fundamental endomorphism.
We define the {\em horizontal curve} $h(t)$ as 
the map that takes $t$ to the kernel of the fundamental projection $P(t)$.

\end{Def}

The horizontal curve is clearly equivariant: if  $T \in Gl(kn)$ then the horizontal curve  of $Tl(t)$ is $Th(t)$.
Observe that since the fundamental endomorphism depends only on the curve on the Grassmannian, the same holds for all its time
derivatives.
In particular, the curve of reflections $D(t)$ and the curve of projections $P(t) := \frac{(I-D(t))}{2}$ depend only
on the curve on the Grassmannian.

We will study now  the second derivative $\ddot{F}$, for this observe that $\dot{P}=-\frac{1}{k}\ddot{F}$.

\begin{Prop}\label{pponto}
Let $\ell(t)$ be a fanning curve and let $h(t)$ be its horizontal curve.
If $P(t)$ is the projection onto $v(t)$ with kernel $h(t)$, then $\dot{P}(t)$ has the following
properties:

1)$\dot{P}(t)$ maps $h(t)$ into $v(t)$.

2)$\dot{P}(t)$ maps $v(t)$ into $h(t)$. 

\end{Prop}

\dem

Differentiating the identity $P(t)^2= P(t)$, we have 
$$\dot{P}(t)P(t)= (I-P(t))\dot{P}(t),$$ where  $I-P(t)$ is the projection onto $h(t)$
with kernel
 $$\spand.$$ Thus the equation
$$0 = \dot{P}(t)P(t)(h(t))= (I- P(t))\dot{P}(t)(h(t))$$
implies that $\dot{P}(t)$ maps $h(t)$ into $\spand$, and this proves the first item.

For the second item, observe that
$$\dot{P}(t)\ell(t))= \dot{P}(t)P(t)\ell(t))= (I- P(t))\dot{P}(t)\ell(t)),$$
which implies that the subspace $\dot{P}(t)\ell(t))$ is contained in $h(t)$. Similarly \\
$\dot{P}(t)(A^{(i)}(t))=  (I- P(t))\dot{P}(t)(A^{(i)}(t))$, for all $i \in \{1,...,k-2\}$ , then $\dot{P}(t)(A^{(i)}(t))$,
for all $i \in \{1,...,k-2\}$, is contained in $h(t)$.

It follows from the proof of Proposition \ref{d} that

$
\begin{array}{l}
1)\dot{F}(t)A(t)=-A(t),\\
2)\dot{F}(t)A(t)^{(i)}=-A(t)^{(i)},$ for all $ i \in \{1,...,k-2\},\\
3)\dot{F}(t)A(t)^{(k-1)}=(k-1)A(t)^{(k-1)}-F(t)A(t)^{(k)}.
\end{array}
$

Differentiating the first equation, we obtain that $ \ddot{F}(t)A(t)=0$, then \\ $\dot{P}(t)A(t)=0$; and similarly we have that  $ \ddot{F}(t)A(t)^{(i)}=0$, 
then \\
$\dot{P}(t)A(t)^{(i)}=0$ for $i \in \{1,2,...,k-3\}$.
Differentiating $\dot{F}(t)A(t)^{(k-2)}=-A(t)^{(k-2)}$, we have that $\ddot{F}(t)A(t)^{(k-2)} + \dot{F}(t)A(t)^{(k-1)}= -A(t)^{(k-1)}$; and using 3)
we obtain  $\ddot{F}(t)A(t)^{(k-2)}=-k A(t)^{(k-1)} + F(t)A(t)^{(k)}$, and consequently
$\dot{P}(t)A(t)^{(k-2)}= A(t)^{(k-1)} - \frac{1}{k}FA^{(k)}$. Since the columns of $FA^{(K)}$ are linear combinations of
those of $A(t)$, $\dot{A}(t)$, ..., $A(t)^{(k-2)}$; and $A(t)$, $\dot{A}(t)$, ...,$A(t)^{(k-2)}$ and $A(t)^{(k-1)}$ are linearly independents, it follows that
$\dot{P}(t)A(t)^{(k-2)}$ has rank $n$. Therefore $\dot{P}(t)A(t)^{(k-2)}$ spans $h(t)$.  \qed

\medskip

\noindent{\em Remark.}
Observe that the proof of proposition \ref{pponto} gives a somewhat more precise information on how 
$\dot{P}(t)$ acts on the associated flags: for any frame, we have the nested flag and 
$\dot{P}(t)$ restricted to each subspace $Span\{A(t), \dot{A}(t), \dots,  A(t)^{(r)}\}$ is zero for each $r<k-2$ and
$\dot{P}(t)$ maps the quotient $v(t)/\{Span\{A(t), \dot{A}(t), \dots,  A(t)^{(k-3)}\}$  isomorphically 
onto $h(t)$; if the frame $A(t)$ is normal, then 
$\dot{P}(t)$ restricted to each subspace $Span\{A(t)\}, Span\{\dot{A}(t)\}, \dots,  Span\{A(t)^{(r)}\}$ 
is zero for each $r<k-2$ and
$\dot{P}(t)$ maps $\{A(t)^{(k-2)}\}$  isomorphically 
onto $h(t)$. 

\subsection{The horizontal derivative } \label{horder}
\begin{Def}
The horizontal derivative of a fanning frame $A(t)$ is the curve of frames defined for:
$$ H(t):= (I- P(t))A(t)^{(k-1)} = \dot{P}(t)A(t)^{(k-2)} = A(t)^{(k-1)} -\frac{1}{k}F(t)A(t)^{(k)}=$$ $$= -\frac{1}{k}\ddot{F}(t)A(t)^{(k-2)},$$
and observe that $H(t)$ is the projection of $A(t)^{(k-1)}$ onto $h(t)$.
\end{Def}

If $A^{(k)} + {k \choose 1}A^{(k-1)}P_1 + {k \choose 2}A^{(k-2)}P_2 +...+{k \choose k-1}\dot{A}P_{k-1}+ AP_{k} = 0$, then
$$ H = A(t)^{(k-1)} + \frac{1}{k}F(t)\left(  {k \choose 1}A^{(k-1)}P_1 +...+{k \choose k-1}\dot{A}P_{k-1}+ AP_{k}\right) ,$$

and using that $F(t)A(t)^{(k-i)}= (k-i)A(t)^{(k-i-1)}$ and ${k\choose i}\frac{k-i}{k}={k-1 \choose i},$
we have
\begin{equation}\label{derhor}
H(t) = A^{(k-1)} + {k-1 \choose 1}A^{(k-2)}P_1 +...+{k-1 \choose k-2}\dot{A}P_{k-2}+ AP_{k-1}.
\end{equation}

\begin{Prop}
The horizontal derivative $H(t)$ of a fanning frame $A(t)$ satisfies the following properties:
\begin{enumerate}
\item  If $X(t)$ is a smooth curve of invertible $n \times n$ matrices, the horizontal derivative of $A(t)X(t)$ is $H(t)X(t)$.
\item  If $T$ is a invertible linear transformation from $\mathbb{R}^{kn}$ to itself, the horizontal derivative of $TA(t)$ is $TH(t)$.
\end{enumerate}
\end{Prop}

\dem

The first property follows from $H(t)=\dot{P}(t)\dot{A}(t)$ e and the proposition \ref{pponto}, since we have $$\dot{P}(t)\frac{d}{dt}(A(t)X(t))= \dot{P}(t)\dot{A}(t)X(t)+\dot{P}(t)A(t)\dot{X}(t)= \dot{P}(t)\dot{A}(t)X(t).$$

The second one is obtained directly from equation \ref{derhor}. \qed

We saw that the curve $\dot{F}(t)$ is a curve of linear transformations with two eigenvalues, $-1$ and $k-1$.  The $-1$-eigenspace  is the vertical space $v(t)$ and the $k-1$-eigenspace is $h(t)$, spanned by the horizontal derivative $H(t)$.  Therefore   $(A(t)|...|A(t)^{(k-2)}|H(t))$ is a natural lift of the curve $\ell(t)$ to
$GL(kn)$, depending on the $k$-jet of the curve. It is worth remarking that once one has a normal form, another 
possible ``natural" lift is given by just using plain derivatives  $(A(t)|...|A(t)^{(k-2)}|A(t)^{(k-1)})$; however one 
loses track of the geometry of the canonical reflection this way. In fact in the real projective plane case, Cartan 
(\cite{cartan})
uses the first lift, that is the one given by the horizontal curve as last columns.

\subsection{The Jacobi Endomorphism} \label{Jacobi}

Taking the derivative of the fundamental reflection, we reach the desired invariant:

\begin{Def}
Let $\ell(t)$ be a fanning curve, $F(t)$ be its fundamental endomorphism and $h(t)$ be the horizontal curve associated to $\ell(t)$.
The Jacobi endomorphism of $\ell(t)$ is defined as $K(t):= \ddot{F}(t)^2/k^2$.
If $P(t) = \frac{(I-D(t))}{2}$, then $K(t)=\dot{P}(t)^2$.
\end{Def}
If $A(t)$ be a fanning frame spanning  $\ell(t)$ and $H(t)$ be its horizontal derivative, we can observe that \begin{equation}\label{jac}
P(t)\dot{H}(t)=-\dot{P}(t)H(t)=-\dot{P}(t)^2A^{(k-2)}(t)=-K(t)A^{(k-2)}(t).
\end{equation}

\begin{The}
Let $\ell(t)$ be a fanning curve in $Gr(n,\mathbb{R}^{kn})$ and let $h(t)$ be its horizontal curve. The Jacobi endomorphism satisfies the following properties:
\begin{enumerate}
\item At each value of $t$, the endomorphism $K(t)$ preserves the decomposition 
 $ \mathbb{R}^ n = v(t) \oplus h(t)$.
\item If $T$ is a transformation in $GL(kn)$, then the Jacobi Endomorphism of $Tl(t)$ is $TK(t)T^{-1}$.
\end{enumerate}
\end{The}

\dem

The first item follows from the proposition \ref{pponto} and from the expression $K(t)=\dot{P}(t)^2$. The second item follows follows from the action of $GL(kn)$ in the fundamental endomorphism. \qed

\vspace{0.2cm}

Observe that, when we consider $A(t)$ a normal frame as in section 2, this frame satisfies:
\begin{equation}\label{normalk}
A^{(k)} + {k \choose 2}A^{(k-2)}\kappa (t)+ {k \choose 3}A^{(k-3)}h_{1}(t) +...+ Ah_{k-2}(t) = 0.
\end{equation}

And in this case we have the horizontal derivative, projecting  $A(t)^{(k-1)}$ onto horizontal curve, takes the form
\begin{equation}\label{normalH}
H(t)= A^{(k-1)} + {k-1 \choose 2}A^{(k-3)}\kappa + {k-1 \choose 3}A^{(k-4)}h_{1} +...+ {k-1 \choose k-1}Ah_{k-3},
\end{equation}
where $\kappa (t) = 1/2\{A(t),t\}.$

\medskip

\begin{The}\label{jacobik}
Let $A(t)$ be a normal frame and $H(t)$ be its horizontal derivative, the matrix of the Jacobi Endomorphism $K(t)$ associated to $A(t)$ in the base of $\mathbb{R}^{kn}$ formed by the columns of $(A(t)|...|A(t)^{(k-2)}|H(t))$ is
$$\left(\begin{array}{cccccc} 0& 0& ...& 0& h_{k-2}-h_{k-3}' & 0 \\ 0&0& ...& 0& {k-1 \choose k-2}(h_{k-3}-h_{k-4}')&0 \\ \vdots& \vdots & &\vdots & \vdots &\vdots\\ 0 & 0 &...&0& {k-1 \choose 2} (h_1 -\kappa')& 0 \\ 0&0&...&0& (k-1)\kappa &0 \\ 0&0&...&0& 0& (k-1)\kappa \end{array} \right),$$
where $\kappa(t)=\frac{1}{2}\{A(t),t\}$.
\end{The}

\dem

First we have that $K(t)A(t)=0, K(t)\dot{A}(t)=0, ..., K(t)A(t)^{(k-3)}=0$, from   proposition \ref{pponto}.
The proof then follows from equations \ref{jac} and \ref{normalH}. Differentiating $H(t)$ in \ref{normalH}, replacing
$$A^{(k)}=-{k \choose 2}A^{(k-2)}\kappa (t)- {k \choose 3}A^{(k-3)}h_{1}(t) -...- Ah_{k-2}(t)$$
and using the property that
$${i-1 \choose j}-{i\choose j} = - {i-1 \choose j-1},$$
we obtain that

\begin{eqnarray}
\dot{H}(t)&=&-{k-1 \choose 1}A^{(k-2)}\kappa - {k-1 \choose 2}A^{(k-3)}(h_1-\kappa')- ... - {}\nonumber \\
& & {}-{k-1 \choose k-2}\dot{A}(h_{k-3}-h'_{k-4})-{k-1 \choose k-1}A (h_{k-2}-h'_{k-3}). \nonumber
\end{eqnarray}

From equation \ref{jac}, we have that $K(t)A(t)^{(k-2)}=-P(t)\dot{H}(t)$, so

\begin{eqnarray}
K(t)A(t)^{(k-2)}&=&{k-1 \choose 1}A^{(k-2)}\kappa + {k-1 \choose 2}A^{(k-3)}(h_1-\kappa')+ ... + {} \nonumber \\
& & {}+{k-1 \choose k-2}\dot{A}(h_{k-3}-h'_{k-4}) +{k-1 \choose k-1}A (h_{k-2}-h'_{k-3}). \nonumber
\end{eqnarray}

Then we just observe that $$K(t)H(t)= \dot{P}(t)(\dot{P}(t)H(t))= \dot{P}(t)(-P(t)\dot{H}(t)), $$
therefore $K(t)H(t)=H(t)(k-1)\kappa(t)$, as claimed. \qed

\medskip

The advantage of taking the square of $\dot{P}$ is that it preserves the vertical-horizontal decomposition; however it might 
be useful to consider $\dot{P}$ itself:

\begin{Prop} \label{Jacobi-in-base}
Let $A(t)$ be a normal frame and $H(t)$ be its horizontal derivative, the matrix of the transformation $\dot{P}(t)$ associated to  $A(t)$ in the base of $\mathbb{R}^{kn}$ formed by the columns of $(A(t)|...|A(t)^{(k-2)}|H(t))$ is
$$\left(\begin{array}{cccccc} 0& 0&...& 0& 0 &h_{k-2}-h_{k-3}'  \\ 0&0& ...& 0& 0&{k-1 \choose k-2}(h_{k-3}-h_{k-4}') \\ \vdots& \vdots & &\vdots & \vdots &\vdots\\ 0 & 0 &...&0& 0& {k-1 \choose 2} (h_1 -\kappa') \\ 0&0&...&0&0 & (k-1)\kappa \\ 0&0&...&0& I& 0 \end{array} \right),$$
where $I$ represents the identity matrix.
\end{Prop}

\dem

From   proposition \ref{pponto} we have that $\dot{P}(t)A(t)=0, ..., \dot{P}(t)A(t)^{(k-3)}=0$ and $\dot{P}(t)A(t)^{(k-2)}=H(t)$, and from the preceding proof, we observe that

\begin{eqnarray}
\dot{P}(t)H(t)&=&{k-1 \choose 1}A^{(k-2)}\kappa + {k-1 \choose 2}A^{(k-3)}(h_1-\kappa')+ ... + {} \nonumber \\
& & {}+{k-1 \choose k-2}\dot{A}(h_{k-3}-h'_{k-4}) +{k-1 \choose k-1}A (h_{k-2}-h'_{k-3}). \nonumber
\end{eqnarray} \qed

\medskip

\obs In the case $k=3$ and $n=1$, that is, the curves in the projective plane $\mathbb{R}P^{2}$, the matrix of theorem \ref{jacobik} is the same that Cartan found in \cite{cartan}. It is interesting to note that, in contrast to the case 
$n=2$ (where there is just a sign difference, see \S 8.3 in \cite{duran}), the matrix of proposition 
\ref{Jacobi-in-base} is {\em not} given 
by pulling back the Maurer-Cartan form by the lift 
$\mathbf{A}(t)=(A(t)|...|A(t)^{(k-2)}|H(t))$ nor by the other plausible lift 
$\tilde{\mathbf{A}}(t)=(A(t)|...|A(t)^{(k-2)}|A(t)^{(k-1)})$. Indeed, for example for $k=4$, we have 
\begin{eqnarray*}
	\mathbf{A}^{-1} \mathbf{A}' &=&\left(\begin{array}{cccc} 0&  0& -h_{1} &  h'_{1}-h_{2}\\ 1& 0& -3\kappa & 3(\kappa'-h_{1})\\  0 & 1 & 0& -3\kappa \\ 0&0& 1& 0  \end{array} \right)\, , \\
	\tilde{\mathbf{A}}^{-1}\tilde{\mathbf{A}}'
	&=&
	\left(\begin{array}{cccc} 0&  0& h_{2} &  3h_{1}.\kappa\\ 0& 0& 3h_{1} & 9\kappa^{2}\\  0 & 0 &  3\kappa& 0 \\ 0&0& 0& 3\kappa  \end{array} \right) \, .
\end{eqnarray*}

\medskip

The last two results explicitly relate the invariants obtained by the fundamental endomorphism and its derivatives and those obtained by the normal forms inspired by the classical invariant theory of projective ODEs. In the next 
section we shall see that the fundamental endomorphism follows naturally and rigidly from an Adjoint representation 
of the space of jets of fanning curves.

\section{Geometry of jets of fanning curves in the Grassmannian}\label{k-jets}

Last section shows that the fundamental endomorphism and related constructions furnish conjugation equivariant maps.
Here we study these maps, especially their uniqueness, as equivariant maps 
 from jets of fanning curves onto the Lie algebra $\lieglkn$. 

First we describe the space $J^{r}_f(\mathbb{R};Gr(n,\mathbb{R}^{kn}))$ of $r$-jets of 
fanning curves 
on the Grassmannian $Gr(n,\mathbb{R}^{kn})$  as 
the quotient of the space $ J^{r}_f(\mathbb{R}; M_{kn\times n})$ of $r$-jets of fanning frames 
by the action of the group $J^{r}(\mathbb{R};Gl(n))$ of $r$-jets of smooth curves of 
invertible $n\times n$ matrices; 
when $k=1$ this is the standard action  $X \cdot A \rightarrow AX$, 
and we extend it to $r$-jets by reppeatly applying Leibnitz's rule. For details of jet groups and their actions, see \cite{Kolar}.

For example, in the case $k=3$, the actions of $J^{2}(\mathbb{R};Gl(n))$ and $ Gl(3n)$  on $ J^{2}_f(\mathbb{R}; M_{3n\times n})$ are given by:
$$( A,\dot{A}, \ddot{A}) \cdot (X,\dot{X},\ddot{X}) = (AX, \dot{A}X +A\dot{X}, \ddot{A}X +2\dot{A}\dot{X}+A\ddot{X})  $$
and
$$T \cdot (A,\dot{A}, \ddot{A}) = (TA,T\dot{A}, T\ddot{A}).$$

In general we have $J^{r}(\mathbb{R};Gl(n))$ and $ Gl(kn)$  acting on $ J^{r}_f(\mathbb{R}; M_{kn\times n})$ 
in the following way:
$$( A,\dot{A},...,A^{(r)}) \cdot (X,\dot{X},..., X^{(r)})=$$ $$ = (AX, \dot{A}X +A\dot{X},..., {r\choose0}A^{(r)}X+ 
{r\choose1}A^{(r-1)}\dot{X}+...+ {r\choose r}AX^{(r)} ) $$
and
$$T \cdot (A,\dot{A},...,A^{(r)}) = (TA,T\dot{A},..., TA^{(r)}).$$

The first use of this description is to show the transitivity of the $GL(kn)$-action on the spaces 
of $k-1$ and $k$-jets of curves in the divisible Grassmannian $Gr(n,kn)$:

\begin{Prop}\label{actstrans}
The group of invertible linear transformations of $\mathbb{R}^{kn}$ acts transitively 
on the space $J^{k-1}_f(\mathbb{R}; M_{kn\times n})$ and, a fortiori, on $J^{k-1}_f(\mathbb{R};Gr(n,\mathbb{R}^{kn})).$

\end{Prop}

\dem

If $ (A| \dot{A}| \ddot{A} | \dots |A^{k-1}) \in J^{k-1}_f(\mathbb{R}; M_{kn\times n})$ then we choose 
\[
T:= (A| \dot{A}| \ddot{A} | \dots |A^{k-1}) \in Gl(kn)
\]

so we have

$$T^{-1} \cdot (A| \dot{A}| \ddot{A} | \dots |A^{k-1}) =
\left( \left(\begin{array}{c} I \\ 0 \\\vdots \\ 0 \end{array} \right),
\left(\begin{array}{c} 0 \\ I \\ \vdots\\ 0 \end{array} \right),
\cdots,
\left(\begin{array}{c} 0 \\ 0 \\ \vdots\\ I \end{array} \right) \right).  $$ \qed

\medskip

Let us look now at the space of $k$-jets of $Gr(n,\mathbb{R}^{kn})$. 
Here the action of $ Gl(kn)$ on the space of $k$-jets of fanning frames is not transitive;
we have that the $k$-jets $(A,\dot{A}, ...,A^{(k)})$ and  $(B, \dot{B},..., B^{(k)})$ are in the 
same $Gl(kn)$-orbit if and only if the matrices $(B|\dot{B}|...|B^{(k-1)})^{-1}B^{(k)}$ and 
$(A|\dot{A}|...|A^{(k-1)})^{-1}A^{(k)}$ are equal. However, we still have that $Gl(kn)$ acts 
transitively on $J^{k}_f(\mathbb{R};Gr(n,\mathbb{R}^{kn}))$:

\begin{Prop}\label{actstrans-k}
The group of invertible linear transformations of $\mathbb{R}^{kn}$ acts 
transitively on the space of $k$-jets of fanning curves in $Gr(n,\mathbb{R}^{kn})$.
\end{Prop}

\dem

All that needs to be shown is  that the joint action of $Gl(kn)$ and $J^{k}(\mathbb{R};Gl(n))$ on the space of $k$-jets of fanning frames is transitive. But if $$(A,\dot{A}, ...,A^{(k)}) \in  J^{k}_f(\mathbb{R}; M_{kn\times n})$$ and $$A^{(k)} + {k \choose 1}A^{(k-1)}P_1 + {k \choose 2}A^{(k-2)}P_2 +...+{k \choose k-1}\dot{A}P_{k-1}+ AP_{k} = 0,$$
then if we act on $(A,\dot{A}, ...,A^{(k)})$, on the left by the matrix
$$(A|\dot{A}+AP_1| A^{(2)}+2\dot{A}P_1 +AP_2|...| \ast )$$
where $\ast =A^{(k-1)} + {k-1 \choose 1}A^{(k-2)}P_1 +...+{k-1 \choose k-2}\dot{A}P_{k-2}+ AP_{k-1})^{-1}$ and on the right by the k-jet $(I,P_1,...,P_{k})$, we get
$$\left( \left(\begin{array}{c} I \\ 0 \\ 0\\ \vdots\\0 \end{array} \right),
\left(\begin{array}{c} 0 \\ I\\0\\\vdots\\ 0 \end{array} \right),...,
\left(\begin{array}{c} 0 \\ 0\\0\\ \vdots \\ I \end{array} \right)
\left(\begin{array}{c} 0 \\ 0\\0\\ \vdots \\ 0 \end{array} \right) \right).  $$ \qed

\subsection{Uniqueness}

In this section we show that the Fundamental Endomorphism, horizontal derivative, etc; are essentially 
unavoidable if we want to represent the $GL(kn)$-action on jets as the Adjoint.

 We begin by
characterizing the fundamental endomorphism for curves in  $Gr(n,\mathbb{R}^{kn})$. 
In order to organize the proof, we need the following lemma whose proof is a matrix computation:

\begin{lem} \label{tecnico}
The matrix
$$ \left(
\begin{array}{cccccc} %
{0\choose0} a_0 I & {1\choose0} a_1 I & {2\choose0} a_2 I & ... & {k-2\choose0} a_{k-2} I & {k-1\choose0} a_{k-1} I  \\
                &                     &                   &     &                         &                          \\
0               & {1\choose1} a_0 I   & {2\choose1} a_1 I & ... & {k-2\choose1} a_{k-3} I & {k-1\choose1} a_{k-2} I  \\
                &                     &                   &     &                         &                          \\
0               & 0                   & {2\choose2} a_0 I & ... & {k-2\choose2} a_{k-4} I  & {k-1\choose2} a_{k-3} I  \\
\vdots          & \vdots              & \vdots            &     & \vdots                  & \vdots \\
0               & 0                   & 0                 & ... & {k-2\choose k-2}a_0 I   & {k-1\choose k-2} a_{1} I \\
                &                     &                   &     &                         &                       \\
0               & 0                   & 0                 & ... & 0                       & {k-1\choose k-1} a_0 I
\end{array} \right)$$
equals  $a_0I + a_1F + \frac{a_2F^2}{2}+...+\frac{a_{k-1}F^{k-1}}{k-1!},$
where
$$
 F= \left(
\begin{array}{cccccc} %
0 & I & 0 & ... &0 &0 \\
0  & 0 & 2I& ...& 0&0 \\
\vdots & \vdots & \vdots & &\vdots & \vdots \\
0 & 0 & 0 & ...& (k-2)I &0 \\
0 & 0 & 0 & ...& 0& (k-1)I\\
0 & 0 & 0 & ...& 0& 0
\end{array} \right).
$$

\end{lem}

\medskip

\begin{The}\label{endfund}
A map $$J^{k-1}_f(\mathbb{R};Gr(n,\mathbb{R}^{kn}))\rightarrow \mathcal{G}l(kn)$$ is equivariant with respect to the $Gl(kn)$ action if and only if it is of the form $$ a_0I + a_1F + \frac{a_2F^2}{2}+...+\frac{a_{k-1}F^{k-1}}{k-1!},$$ where $I$ is the identity matrix, $a_i$, for all $i$, are real numbers, and $$ F= \mathbf{A}(t)  \left(
\begin{array}{cccccc}
0 & I & 0 & ... &0 &0 \\
0  & 0 & 2I& ...& 0&0 \\
\vdots & \vdots & \vdots & &\vdots & \vdots \\
0 & 0 & 0 & ...& (k-2)I &0 \\
0 & 0 & 0 & ...& 0& (k-1)I\\
0 & 0 & 0 & ...& 0& 0
\end{array} \right) \mathbf{A}(t)^{-1},$$ with $\mathbf{A}(t) = \jusk$.

\end{The}

\dem

The proof is divided in two parts: first, that in the right basis, the matrix representing the map has to be constant. Then, 
we show that the entries of this matrix are the correct ones to give the desired result.

\noindent{\em First part:}
Let $G: J^{k-1}_f(\mathbb{R}; M_{kn\times n})\rightarrow \mathcal{G}l(kn)$ be a map invariant under the action of $J^{k-1}(\mathbb{R};Gl(n))$ and equivariant with respect to the action of $Gl(kn)$.
Writing  $G(A,\dot{A},...,A^{(k-1)})$ in the canonical basis, we obtain
$$  (A|\dot{A}|...|A^{(k-1)}) \left(
\begin{array}{c}
G_{ij}(A,\dot{A},...,A^{(k-1)}) \end{array} \right)_{k\times k} (A|\dot{A}|...|A^{(k-1)})^{-1},$$
where $G_{ij}$ are blocks $n \times n$.

The equivariance implies that
\[
G(TA,T\dot{A},...,TA^{(k-1)}))= TG(A,\dot{A},...,A^{(k-1)})T^{-1},
\] then, $\forall T \in Gl(kn)$,
$$\left(
\begin{array}{c}
G_{ij}(A,\dot{A},...,A^{(k-1)}) \end{array} \right)_{k\times k} = \left(
\begin{array}{c}
G_{ij}(TA,T\dot{A},...,TA^{(k-1)}) \end{array} \right)_{k\times k}.$$

Since $Gl(kn)$ acts transitively on $J^{k-1}_f(\mathbb{R}; M_{kn\times n})$, then the blocks $G_{ij}$ of $n \times n$ matrices are constant.

\noindent{\em Second part:}
By induction on $k$, assume that the matrices $G_{ij}$ of the map  $$J^{k-2}_f(\mathbb{R};Gr(n,\mathbb{R}^{(k-1)n}))\rightarrow \mathcal{G}l((k-1)n)$$ satisfies the Lemma. Now, using the invariance under the action of $J^{k-1}(\mathbb{R};Gl(n))$, we need to conclude that the matrices $G_{ij}$  of the map  $J^{k-1}_f(\mathbb{R};Gr(n,\mathbb{R}^{kn}))\rightarrow \mathcal{G}l(kn)$ has the form of the Lemma.

The invariance under the action of $J^{k-1}(\mathbb{R};GL(n))$ implies, for all  $X(t) \in Gl(n) $, that:

\begin{eqnarray}\label{prim}
\end{eqnarray}
$$\left(
\begin{array}{cccccc}
{0\choose0} X  & {1\choose0} \dot{X} & {2\choose0} \ddot{X}& ... & {k-2\choose0} X^{(k-2)}& {k-1\choose0} X^{(k-1)}  \\
               &                     &                     &     &                        &                       \\
0              & {1\choose1} X       & {2\choose1} \dot{X} &...  & {k-2\choose1} X^{(k-3)}&{k-1\choose1} X^{(k-2)}  \\
\vdots         & \vdots              & \vdots              &     & \vdots                 & \vdots \\
0              & 0                   & 0                   & ... & {k-2\choose k-2}X      & {k-1\choose k-2}\dot{X} \\
               &                     &                     &     &                        &                 \\
0               & 0                  & 0                  & ... & 0                      & {k-1\choose k-1} X
\end{array} \right)  \left(
\begin{array}{c}
G_{ij} \end{array} \right)_{k\times k}= $$
$$= \left(
\begin{array}{c}
G_{ij} \end{array} \right)_{k\times k} \left(
\begin{array}{cccccc}
{0\choose0} X  & {1\choose0} \dot{X} & {2\choose0} \ddot{X}& ... & {k-2\choose0} X^{(k-2)}& {k-1\choose0} X^{(k-1)}  \\
               &                     &                     &     &                        &                       \\
0              & {1\choose1} X       & {2\choose1} \dot{X} &...  & {k-2\choose1} X^{(k-3)}&{k-1\choose1} X^{(k-2)}  \\
\vdots         & \vdots              & \vdots              &     & \vdots                 & \vdots \\
0              & 0                   & 0                   & ... & {k-2\choose k-2}X      & {k-1\choose k-2}\dot{X} \\
               &                     &                     &     &                        &                 \\
0               & 0                   & 0                  & ... & 0                      & {k-1\choose k-1} X
\end{array} \right).$$

\bigskip

Looking at the last line of these products, we obtain the relations:

$XG_{k1} = G_{k1}X $

$XG_{k2} = G_{k1} \dot{X} + G_{k2}X $

$\hspace{1cm} \vdots $

$XG_{k,k-1} = {k-2\choose0}G_{k1} X^{(k-2)}+ \dots + {k-2\choose k-2}G_{k,k-1} X.$

Then $G_{k1}, G_{k2}, ..., G_{k,k-1}$ must be zero.

Therefore $G$ has the form $\left(
\begin{array}{ccccc} %
 & &                    &  G_{1k} \\
 & \ast &               &  G_{2k}  \\
 & &                    & \vdots \\
0 & ...   &  0             &  G_{kk}
\end{array} \right),$
where $\ast$ just depends of the first $k-1$ lines and columns of the matrix $(G_{ij})_{k\times k}$, and of $J^{k-2}(\mathbb{R};Gl(n)).$

Using the hypotheses, $G$ has the form:

\begin{eqnarray} \label{seg}
\left(
\begin{array}{cccccc} %
{0\choose0} a_0 I & {1\choose0} a_1 I & {2\choose0} a_2 I & ... & {k-2\choose0} a_{k-2} I & G_{1k}  \\
                  &                   &                   &     &                         &         \\
0               & {1\choose1} a_0 I   & {2\choose1} a_1 I & ... & {k-2\choose1} a_{k-3} I & G_{2k}  \\
                &                     &                   &     &                         &         \\
0               & 0                   & {2\choose2} a_0 I & ... & {k-2\choose2} a_{k-4} I   & G_{3k}  \\
\vdots          & \vdots              & \vdots            &     & \vdots                  & \vdots \\
0               & 0                   & 0                 & ... & {k-2\choose k-2}a_0 I   & G_{k-1,k} \\
                &                     &                   &     &                         &        \\
0               & 0                   & 0                 & ... & 0                       & G_{kk}
\end{array} \right)
\end{eqnarray}

where, using \ref{prim} and \ref{seg}, $G_{ik}$ satisfies:

$$
{0\choose0}XG_{1k}+ {1\choose0} \dot{X} G_{2k} +\cdots + {k-1\choose0}X^{(k-1)} G_{kk}  = $$ $$= {0\choose0}a_0{k-1\choose0} X^{(k-1)}+  \cdots + {k-2\choose0}a_{k-2}{k-1\choose k-2} \dot{X}+ G_{1k} {k-1\choose k-1}  X $$

\centerline{\vdots}

$${k-3\choose k-3} XG_{k-2,k}+ {k-2\choose k-3} \dot{X} G_{k-1,k} + {k-1\choose k-3} X^{(2)} G_{kk} = $$ $$= {k-3\choose k-3}a_0{k-1\choose k-3} X^{(2)}+ {k-2\choose k-3}a_1{k-1\choose k-2}\dot{X}  + G_{k-2,k} {k-1\choose k-1}  X  $$

$${k-2\choose k-2} XG_{k-1,k}+ {k-1\choose k-2} \dot{X} G_{kk} = {k-2\choose k-2}a_0{k-1\choose k-2}\dot{X} + G_{k-1,k} {k-1\choose k-1}  X $$

$${k-1\choose k-1} X G_{kk}  =   G_{kk} {k-1\choose k-1}X $$

Therefore the matrix $\left(
\begin{array}{c}
G_{ij} \end{array} \right)_{k\times k} $  has the form of the Lemma \ref{tecnico}. \qed

\vspace{1.0cm}

Now we characterize the horizontal derivative:

\begin{The}\label{kjethor}
The assignment that sends a fanning curve $\ell(t)$ to its horizontal curve $h(t)$ is characterized by the following four properties:
\begin{enumerate}
\item At each $t$, the subspace $h(t)$ is transversal to $Span\{A(t),...,A^{(k-2)}\}$.
\item The subspace $h(\tau)$ depends only on the $k$-jet of the curve $\ell(t)$ at $t=\tau$.
\item If $T\in Gl(kn)$, the horizontal curve of $Tl(t)$ is $Th(t)$.
\item If $\ell(t)$ is spanned by a curve $A_0+tA_1+...+t^{k-1}A_{k-1}$ in the space of frames $h(t)$, $h(t)$ is constant.
\end{enumerate}
\end{The}

The next lemma, necessary for the proof of  theorem \ref{kjethor},   is the analogue of  lemma 7.5 of (\cite{duran}) for the generalized  horizontal derivative.

\begin{lem}
 A $Gl(kn)$-equivariant map $j:J^{k}_f(\mathbb{R};Gr(n,\mathbb{R}^{kn})) \rightarrow Gr(n,\mathbb{R}^{kn})$ is such that the subspaces $j([(A,\dot{A},...,A^{k})])$ and $Span\{A(t),...,A^{(k-2)}\}$ are always transversal if and only if it is of the form $$[(A,\dot{A},...,A^{(k)})]\longmapsto [H+ c_{k-1}A +c_{k-2}\dot{A}+...+ c_{1}A^{k-2}], $$
where $c_1,...,c_{k-1}$ are real numbers and $H $ is the horizontal derivative defined in \ref{derhor}.
\end{lem}

\dem

Since $(A,\dot{A},...,A{(k)})\mapsto H$ is the horizontal derivative, for any real numbers $c_1,...,c_{k-1}$, the subspace $[H+ c_{k-1}A +c_{k-2}\dot{A}+...+ c_{1}A^{k-2}]$ is transversal to  $Span\{A(t),...,A^{(k-2)}\}$. And the $Gl(kn)$-equivariance follows from the properties of the horizontal derivative.

Conversely, let us define $P: J^{k}_f(\mathbb{R}; M_{kn\times n})\rightarrow \mathcal{G}l(kn)$ as the map whose value at a $k$-jet $(A,\dot{A},...,A^{(k)})$ is the projection with range $Span\{A(t),...,A^{(k-2)}\}$ and kernel $J([(A,\dot{A},...,A^{(k)})])$. The map $P$ has the following properties:
\begin{enumerate}
 \item $P(A,...,A^{(k)})^2 = P(A,...,A^{(k)})$;
 \item $P(A,...,A^{(k)})A^{(i)}=A^{(i)}$, for $i=1,...,k-2$;
 \item $P(TA,...,TA^{(k)})=TP(A,...,A^{(k)})T^{-1}$;
 \item $P((A,\dot{A},...,A^{(k)}) \cdot (X,\dot{X},...,X^{(k)}))=P(A,\dot{A}, ...,A^{(k)})$.
\end{enumerate}

Using (1) and (2), we have that there exists $k-1$ functions 

\centerline{$R_1(A,\dot{A},...,A^{(k)})$, $R_2(A,\dot{A},...,A^{(k)})$, ..., $R_{k-1}(A,\dot{A},...,A^{(k)})$}
 with values in the space of $n\times n$ matrices such that $P(A,...,A^{(k)})$ is equal to
$$( A(t)|\dot{A}(t)|...|A^{(k-1)}(t)) \left(
\begin{array}{ccccc} %
I   & 0  & ... & 0 & R_1  \\
0   & I  & ... & 0 & R_2   \\
\vdots & \vdots  &  & \vdots & \vdots \\
0  & 0  & ... & I   & R_{k-1} \\
0  & 0  & ... & 0   & 0
\end{array} \right) ( A(t)|\dot{A}(t)|...|A^{(k-1)}(t))^{-1}.$$

Since $(A|\dot{A}|...|A^{(k-1)})^{-1}A^{(k)}$ is the complete invariant for the action of $Gl(kn)$ on the $k$-jets of fanning frames (Proposition \ref{actstrans}), and   property
(3) implies that $R_i(A,\dot{A},...,A^{(k)})$, for all $i\in \{1,...,k-1\}$, depends only of the $Gl(kn)$-orbit; then  $R_i(A,\dot{A},...,A^{(k)})$
depend only of $(A|\dot{A}|...|A^{(k-1)})^{-1}A^{(k)}$.

Moreover,  property (4) is equivalent to the $R_i$, for all $i$, having the following expressions:
$$R_1 = -P_{k-1}-c_{k-1}I,$$ $$R_2= - {k-1 \choose k-2} P_{k-2} - c_{k-2}I,$$ $$\vdots$$ $$R_{k-1} = - {k-1 \choose 1}P_1-c_1I$$

where $P_i$, for all $i$, comes from $H(t)$.
 Therefore, since $j([(A,\dot{A},...,A^{k})])$ is the kernel of $P(A,...,A^{(k)})$, it must be $[H+ c_{k-1}A +c_{k-2}\dot{A}+...+ c_{1}A^{k-2}]$. \qed

\vspace{1.0cm}

\textit{Proof of Theorem \ref{kjethor}}

If a map $J^{k}_f(\mathbb{R};Gr(n,\mathbb{R}^{kn})) \rightarrow Gr(n,\mathbb{R}^{kn})$ is equivariant and it has the property that its image is always transversal to $Span\{A(t),...,A^{(k-2)}\}$ then from the previous lemma
$h(t)= [H+ c_{k-1}A+c_{k-2}\dot{A}+...+c_{1}A^{(k-2)}]$, for any choice of frame $A(t)$ spanning $\ell(t)$.
Since $Gl(kn)$ acts transitively on $J^{k}_f(\mathbb{R};Gr(n,\mathbb{R}^{kn}))$ and the map is  equivariant, 
we just need analyze the map in one point then we determine it. 
When $A(t)$ has the form $A_0+tA_1+...+t^{k-1}A_{k-1}$ in the space of frames, 
then $h(t)$ is $$[A_{k-1}+ c_{k-1}(A_0+tA_1+...+t^{k-1}A_{k-1})+...+ c_{1}((k-2)!A_{k-2}+(k-1)!\;tA_{k-1})]$$ 
and this is constant if and only if $c_1,c_2, ...,c_{k-1}$ are zero. So, as claimed, $h(t)=[H(t)]. $ \qed

\subsection{$k$-jets and Adjoint orbits}

Let us examine more closely the invariants of section \ref{fundamentaletc} as equivariant maps from the space of $r$-jets of curves onto the Lie algebra $\lieglkn$
endowed with the Adjoint action. We will do the interesting case $r=k-1$, where 
we shall see that the fundamental endomorphism is actually an equivariant 
{\em embedding}, thus modelling $J^{k-1}_f(\mathbb{R};Gr(n,\mathbb{R}^{kn}))$  as 
an adjoint orbit; $r=k$, where the $GL(kn)$ is still transitive, and 
$r=k+1$, the first stage where the action ceases to be transitive and we 
shall see how the invariants parametrize the space of orbits. 

\begin{Prop}
 The fundamental endomorphism  
 $$F:J^{k-1}_f(\mathbb{R};Gr(n,\mathbb{R}^{kn}))\rightarrow \lieglkn$$
of theorem \ref{endfund} is a diffeomorphism onto its image, in fact, $F$ is a equivariant embedding of $k-1$-jets  
as an Adjoint orbit Lie algebra of $\lieglkn$.
\end{Prop}

\dem

Since $GL(n,kn)$ acts transitively on $J^{k-1}_f(\mathbb{R};Gr(n,\mathbb{R}^{kn}))$, and the map is equivariant, 
its image is contained in an Adjoint orbit. All we need to check is that given 
$s\in J^{k-1}_f(\mathbb{R};Gr(n,\mathbb{R}^{kn}))$, the isotropy of $s$ and the isotropy of $F(s)$ coincide, 
which holds since both isotropies are 
composed of matrices with the form
$$\left(
\begin{array}{cccccc}
X_1 & X_2 & X_3 & ... & X_{k-1} & X_k \\
0  & X_1 & 2X_2& ...& {k-1 \choose 2} X_{k-2} & {k \choose 2} X_{k-1}   \\
\vdots & \vdots & \vdots & &\vdots & \vdots \\
0 & 0 & 0 & ...& {k-1 \choose k-2}X_2 & {k \choose k-2} X_3 \\
0 & 0 & 0 & ...& X_1 & {k \choose k-1 }X_2\\
0 & 0 & 0 & ...& 0 & X_1
\end{array} \right).$$ \qed

\medskip

\noindent{\bf Remark.} The previous result holds by taking as equivariant map any map of the form 
$ a_0I + a_1F + \frac{a_2F^2}{2}+...+\frac{a_{k-1}F^{k-1}}{k-1!}$, as long as $a_1 \neq 0$. 

We now study the first and second prolongation of the fundamental endomorphism. As soon as 
$r\geq k$, we have the advantage of having enough information 
to define normal frames: 
Let $N^{r}_f(\mathbb{R}; M_{kn\times n}) \subset J^{r}_f(\mathbb{R}; M_{kn\times n})$ be the space 
of $r$-jet of {\em normal} fanning frames. By the results in section \ref{normal}, we have that  $J^{r}_f(\mathbb{R};Gr(n,\mathbb{R}^{kn}))$ is the
quotient of $N_f^{r}(\mathbb{R}; M_{kn\times n})$ 
by the group $GL(n)$ of (constant) $n\times n$ invertible matrices. 

The fundamental projection $P(t)$ gives a map $P:J^{k}_f(\mathbb{R};Gr(n,\mathbb{R}^{kn}))$ into the space of 
projections; more precisely, into the connected component $\Pi(n,kn)$ of the space of projections of $\mathbb{R}^{kn}$ indexed by 
$n = \dim \ker P(t) = \dim h(t) = n$. Recall that everything is linear (as opposed to Euclidean), and 
the map 
\begin{eqnarray*}
	\Pi(n,kn)  &\to& Gr(n,kn)\\
	\pi        &\mapsto& \ker(\pi)
\end{eqnarray*}
is a submersion. The space of linear projections is used, for example, as the classifying space in the category 
of vector bundles endowed with linear connections (\cite{porta-recht}). 

Since the action of $GL(kn)$ is transitive both on $k$-jets of curves and 
the space $\Pi(n,kn)$,  the fundamental projection gives 
a {\em surjective} equivariant map $P:J^{k}_f(\mathbb{R};Gr(n,\mathbb{R}^{kn})) \to \Pi(n,kn)$. This map can be factored through the flags appearing previously in the paper; denoting by $\mathcal{F}(n,kn)$ (resp $\mathcal{D}(n,kn)$) the flag spaces 
of linear chains of subspaces (resp. decompositions) of $\mathbb{R}^{kn}$ of the 
appropriate dimensions, we have the submersions
\[
J^{k}_f(\mathbb{R};Gr(n,\mathbb{R}^{kn}))\stackrel{d}{\to}
\mathcal{D}(n,kn)\to
\mathcal{F}(n,kn)\to
\Pi(n,kn)\to 
Gr(n,kn)\, .
\]

All the arrows with the possible exception of the first one are well understood. In order to grasp the first map $d$, by equivariance and transitivity it is a homogeneous submersion whose typical  fiber is the quotient of the isotropies
$I_{d(x)}/I_x$. The isotropy of a given decomposition in $\mathcal{D}(n,kn)$ is the set of linear transformations that preserve each space, i.e., the $k$-fold product
$\mathcal GL(n)^k$. 

Now let $T\in  GL(kn)$ fixing a $k$-jet $j_k  \in J^{k}_f(\mathbb{R};Gr(n,\mathbb{R}^{kn}))$ of the form given in the transitivity 
proposition \ref{actstrans-k}. By lifting $j_k$ to a normal 
$k$-jet of frames $\mathbb{A}_k\in N^{r}_f(\mathbb{R}; M_{kn\times n}) $ we have that $\mathbb{B}=T\mathbb{A}$ is also a normal $k$-jet, such that 
$Span(A^{(r)}) = Span(B^{(r)}), 0\leq r \leq k-2$ and $Span(H_{\mathbb{A}}) = Span(H_{\mathbb{B}})$. Therefore, there exists a constant, invertible 
$B$ such that $A^{(r)}X = B^{(r)} $ for all $0\leq r \leq k-2$ and also $H_{\mathbb{A}}X= H_{\mathbb{B}}$. That means that $T$ must be a block-diagonal 
matrix $(X,\dots,X)\in \Delta \subset GL(n)^k \subset GL(kn)$. Thus the fiber 
of the map $d$ is the homogenous space $GL(n)^k/\Delta$, which is diffeomorphic
as a differentiable manifold 
to $GL(n)^{k-1}$, a diffeomorphism being realized by the ``''homogenous coordinates"
\[
(g_1,\dots g_k) \mapsto (g_k^{-1}g_1,g_k^{-1}g_2, \dots g_k^{-1}g_{k-1}) \, .
\]

Note that $GL(n)$ sits inside of the set $M_n$ of all $n\times n$ matrices, 
and $GL(n)$ acts diagonally on $M_n^k$. The quotient is a ``non-commutative
projective space" (it is actually $\mathbb{R}P^{k-1}$ when $n=1$) and the fiber 
is the open set that is the intersection of the domains of the homogeneous 
coordinate charts. 

\medskip

Let us finally deal with the space of $k+1$-jets. Here the $GL(kn)$-action is no 
longer transitive and we want to coordinatize the space of orbits. We still 
use the restrict ourselves to normal frames, but additionally, we work on a {\em section}, that is an appropriate submanifold 
of $J^{k+1}_f(\mathbb{R};Gr(n,\mathbb{R}^{kn}))$ that intersects all orbits. 
This is done in order to 
avoid the ambiguity of a choice of basis that translates to conjugation in theorem \ref{mainCongruence}.
Denote by $\{\vec{e}_1,\dots, \vec{e}_{kn}\}$  the canonical basis of $\mathbb{R}^{kn}$. 

\medskip
\begin{Def}
	A  $r$-jet of curves in the divisible Grassmannian 
	to  {\em standard} if its projection to $0$-jets is the plane 
	$Span\{\vec{e}_1, \dots \vec{e}_n\}$. 
	
	A $r$-jet of frames is  {\em standard} if it is normal and its projection 
	to $0$-jets is the frame $(\vec{e}_1, \dots \vec{e}_n)$. 		
\end{Def}

Let us observe that the concept of  standard curves makes sense for all $r$-jets, whereas 
for a standard frame we need $r\geq k$ in order to define normalcy. 
Let us denote by $\mathcal{S}(n,kn)$ (resp.  $\widetilde{\mathcal{S}}(n,kn)$ ) the space of $k+1$-jets of 
standard 
curves in $Gr(n,kn)$ (resp . $k+1$-jet of
curves of  standard frames).

Since normal frames are unique given an initial frame and the initial frame is 
fixed for standard frames, we have

\begin{Prop}\label{standard-mesma-coisa}
	The projection ({\it frame $A$}) $\mapsto$ {\it span(A)}  induces a diffeomorphism $\widetilde{\mathcal{S}}(n,kn)\to \mathcal{S}(n,kn)$.
\end{Prop}

The group $G_0\subset GL(kn)$ that preserves 
 $\mathcal{S}(n,kn)$ is the group of block-upper triangular matrices of the form
 \[
\begin{pmatrix}
  X   & Y  \\
  0   & Z \\
 \end{pmatrix}
  \]
where each $X \in GL(n), Z\in GL(n(k-1)) $. The action of $G_0$ on $\mathcal{S}(n,kn)$ lifts to standard frames as follows:
\[
\begin{pmatrix} %
X   & Y  \\
0   & Z 
\end{pmatrix}
\bullet 
A
=
\begin{pmatrix} %
X   & Y  \\
0   & Z 
\end{pmatrix}
A
X
\] 

It is clear that  $\mathcal{S}(n,kn)$ is indeed a section. Therefore the inclusion 
$\mathcal{S}(n,kn) \hookrightarrow J^{k+1}_f(\mathbb{R};Gr(n,\mathbb{R}^{kn})) $ 
induces a homeomorphism 
$$\mathcal{S}(n,kn)/G_0 \hookrightarrow J^{k+1}_f(\mathbb{R};Gr(n,\mathbb{R}^{kn}))/GL(kn),$$
and by proposition \ref{standard-mesma-coisa}, also a homeomorphism 
$$\widetilde{\mathcal{S}}(n,kn)/G_0 \hookrightarrow J^{k+1}_f(\mathbb{R};Gr(n,\mathbb{R}^{kn}))/GL(kn).$$

We have

\begin{The}
	The map $Q: S(n,kn) \to \lieglkn$ given by the entries of the matrix of theorem 
	\ref{jacobik} induces a homeomorphism between its image and the space of orbits  $J^{k+1}_f(\mathbb{R};Gr(n,\mathbb{R}^{kn}))/GL(kn)$.
\end{The}

\dem

Theorem \ref{mainCongruence} says that two curves are congruent if and only if the respective Schwarzians and matrices $h_j$ are of normal frames lifting them are conjugate  by a constant $n\times n$ matrix $X$. If both curves are standard, then 
$X$ must be the identity. The only missing piece is to substitute ``$k+1$-jet" in place 
of ``curves" in the begginning of the proof; it is not clear at first glance that the 
Schwarzian and the $h_j$ depend on the $k+1$-jet of a curve. But this follows from 
the presentation of the Jacobi endomorphism of Theorem \ref{jacobik} and Proposition \ref{Jacobi-in-base}: the Jacobi endomorphism and its associated matrix in the basis given by $(A,\dot{A}, \dots , H)$ can be computed with using {\em at most} $k+1$ derivatives,  and  one needs {\em at least} $k+1$ derivatives since 
otherwise the Jacobi matrices of \ref{jacobik} and \ref{Jacobi-in-base} would be constant by the transitivity of the action on $r$-jets, $r\leq k$.  \qed

\end{document}